\documentclass[11pt]{article}

\usepackage[utf8]{inputenc}
\usepackage{amsmath,amsthm,amssymb,enumerate,framed,mdwlist}
\usepackage{color,colortab}
\usepackage[square,numbers,sort&compress]{natbib}

\newcommand{\T}{\mathbb T}

\newcommand{\N}{\mathbb{N}}
\newcommand{\Z}{\mathbb{Z}}
\newcommand{\R}{\mathbb{R}}
\newcommand{\conv}{\operatorname{conv}}

\renewcommand{\epsilon}{\varepsilon}

\newtheoremstyle{mythmstyle}
	{\topsep}
	{\topsep}
	{\itshape}
	{}
	{\scshape}
	{.}
	{3pt}
	{}
\theoremstyle{mythmstyle}

\newtheorem{nn}{}[section]
\newtheorem{lemma}[nn]{Lemma}
\newtheorem{theorem}[nn]{Theorem}

\newtheorem{prop}[nn]{Proposition}
\newtheorem{definition}[nn]{Definition}
\newtheorem{claim}[nn]{Claim}

\newtheorem{example}[nn]{Example}

\newtheorem{REMARK}[nn]{Remark}
\newenvironment{remark}{\begin{REMARK}}{\end{REMARK}}

\newenvironment{cpf}{\begin{trivlist} \item[] {\em Proof of Claim.}}{\hspace*{\stretch{1}} $\diamond$ \end{trivlist}}

\def\ve#1{\mathchoice{\mbox{\boldmath$\displaystyle\bf#1$}}
{\mbox{\boldmath$\textstyle\bf#1$}}
{\mbox{\boldmath$\scriptstyle\bf#1$}}
{\mbox{\boldmath$\scriptscriptstyle\bf#1$}}}

\newcommand{\x}{{\ve x}}
\newcommand{\y}{{\ve y}}
\newcommand{\z}{{\ve z}}

\renewcommand{\a}{{\ve a}}

\newcommand{\0}{{\ve 0}}

\newcommand{\bve}{{\ve b}}

\newcommand{\GMI}{\operatorname{GMI}}
\newcommand{\GOM}{\operatorname{GOM}}

\numberwithin{equation}{section}

\allowdisplaybreaks

\begin{document}

\title{Optimal cutting planes from the group relaxations}

\author{Amitabh Basu\footnote{Department of Applied Mathematics and Statistics, The Johns Hopkins University. Supported by the NSF grant CMMI1452820.} \and Michele Conforti\footnote{Dipartimento di Matematica ``Tullio Levi-Civita'', Universit\`a degli Studi di Padova, Italy. Supported by the grants ``SID 2016'' and ``PRIN 2016''.}\and Marco Di Summa\footnotemark[2]\and Giacomo Zambelli\footnote{Department of Mathematics, London School of Economics and Political Science.}}


\maketitle

\begin{abstract}
We study quantitative criteria for evaluating the strength of valid inequalities for Gomory and Johnson's finite and infinite group models  and we describe the valid inequalities that are optimal for these criteria. We justify and focus on the criterion of maximizing the volume of the nonnegative orthant cut off by a valid inequality.

For the finite group model of prime order, we show that the unique maximizer is an automorphism of the {\em Gomory Mixed-Integer (GMI) cut} for a possibly {\em different} finite group problem of the same order.

We extend the notion of volume of a simplex to the infinite dimensional case. This is used to show that in the infinite group model,  the GMI cut maximizes the volume of the nonnegative orthant cut off by an inequality.\end{abstract}

\section{Introduction}


 Cutting planes are important tools to solve integer programming (IP) models. While this technology has seen several revivals and intense research activity since its introduction by Gomory~\cite{MR0182454,gomory1960algorithm,Gomory63,MR0102437}, some basic aspects are not well understood. In particular, deciding {\em which} cutting plane family will be most effective in a particular IP instance has been a thorny problem for the community. Some families, like the {\em Gomory Mixed-Integer (GMI) cuts}, have been enormously useful across all kinds of IP instances~\cite{balas96gomory,balas2008optimizing,dash2010mir,Bixby}, but such empirical observations have never been explained rigorously, to the best of our knowledge.

This problem, known as {\em cut selection}, has become even more confounded by the advent of the theory of so-called {\em cut generating functions}, which is a modern perspective on the theories of Gomory and Johnson~\cite{infinite,infinite2}, and Balas~\cite{bal} from the 1970s. This is because recent work in this area has opened the doors to infinitely many distinct families of cutting planes that are computationally accessible. In a sense, this is great news because we now have not only many more choices, but potentially more powerful cutting planes than used before. On the other hand, this makes the problem of cut selection even more difficult than what it was.

The goal of this paper is to make  progress towards establishing rigorous criteria for evaluating the efficacy of cutting planes, and understanding the structure of the optimal cutting planes under these criteria. For this purpose, we restrict ourselves to the family of cutting planes studied under the various group relaxations. The group relaxations were first introduced by Gomory~\cite{gom}, and then generalized and investigated by Gomory and Johnson~\cite{infinite,infinite2,johnson,tspace}, as a means to provide a unifying framework for deriving cutting planes for mixed-integer optimization problems.

Our investigations seem to provide analytical evidence supporting the computational success of the GMI cuts (see Theorems~\ref{thm:finite-volume} and~\ref{thm:infinite-volume}).
\medskip

 Our presentation  uses notions and  terminology that may be unfamiliar in the integer programming community. However, our approach has a three-fold motivation: 1) this new language enables us to state results in a unified manner and to borrow necessary mathematical machinery from the areas of pure algebra and analysis; 2) the abstraction allows us to present the proofs in a cleaner and more elegant fashion; 3) the added generalization could be useful in the future to build further bridges between algebra, analysis and mixed-integer optimization. We next review this material and then introduce Gomory's finite and infinite group models.

\paragraph{Group theoretic preliminaries.} We recall that a {\em group} is a set $G$ endowed with a binary operation mapping $G\times G$ to $G$, denoted by $+ : G\times G \to G$, which satisfies three properties: 1) $\x + (\y + \z) = (\x + \y) + \z$ for all $\x,\y, \z \in G$ (associativity of $+$), 2) there exists an element $\0 \in G$ such that $\x+\0=\0 + \x  = \x$ for all $\x \in G$ (existence of identity element), and 3) for every $\x	\in G$, there exists an inverse $-\x$ such that $\x + (-\x) = (-\x)+\x= \0$. A group is said to be {\em abelian} if $+$ is commutative, i.e., $\x + \y = \y+\x$ for all $\x,\y \in G$. All groups considered in this paper will be abelian, so we drop this qualification in the remainder.
We define $\x - \y:=\x +(-\y)$ for every $\x,\y\in G$.
Moreover, we will use the notation $k\x$ to denote $\x$ added to itself $k$ times, for any $k \in \N$ and $\x \in G$.

A group $G$ is said to be a {\em topological group} if the set $G$ is also endowed with a topology such that the maps $+: G\times G \to G$ and $\operatorname{inv}:G\to G$, $\x\mapsto-\x$, are continuous with respect to this topology, where $G\times G$ is endowed with the product topology. A topological group is said to be {\em compact} (resp., {\em connected}) if $G$ is a compact (resp., connected) space under its given topology. We will assume that all topological spaces in this paper are Hausdorff, i.e., for any two distinct points there exist disjoint neighborhoods containing these two points.\footnote{For a reference on these terminologies, see~\cite{royden1968real}.}

For any compact, topological group $G$, there exists a unique measure $\mu$ defined on the Borel sets of $G$ such that $\mu(G) = 1$ and $\mu(\x+A)=\mu(A)$ for every Borel subset $A\subseteq G$ (where $\x+A:=\{\x+\a:\a\in A\}$). This measure is called the {\em Haar probability measure}, or {\em Haar measure} for short. It can be also be shown that $\mu(S) = \mu(-S)$ for all subsets $S \subseteq G$. See~\cite[Chapter 22]{royden1968real} for a detailed discussion.

$\T^n := \R^n/\Z^n$ will denote the $n$-dimensional torus. We will use bold-faced letters like $\x,\y,\bve$ etc., to denote elements of a generic group. When the group is a ``one-dimensional" group like $\T^1$ or a finite cyclic group $\Z/q\Z$, we will not use the bold font for elements in the group, i.e., refer to elements using $x,y,b$ etc.


\paragraph{Gomory's group relaxation.} Consider a group $G$ and a fixed $\bve \in G\setminus\{\0\}$. Gomory's group relaxation is defined as the set $I_\bve(G)$ of {\em finite support} functions $y: G \to \Z_+$ (i.e., $y$ takes value $0$ on all but a finite subset of $G$) such that
$$\sum_{\x \in G} y(\x)\x = \bve.$$
The usual cases are with $G = \T^n$ (the $n$-dimensional torus), where it is called the {\em $n$-dimensional infinite group relaxation}, and $G = \Z^n/\Lambda$ where $\Lambda$ is a sub-lattice of $\Z^n$, where it is called an {\em $n$-dimensional finite group relaxation}. These two cases and their use in deriving cutting planes for integer programming have been studied extensively; we refer the reader to the surveys~\cite{basu2016light,basu2016light2,basu2015geometric}.

Note that, as $\bve\ne\0$, the function $y_0$ that takes value 0 on all the elements of $G$ is not in $I_\bve(G)$. In the context of integer programming, the function $y_0=0$ corresponds to the optimal solution of the continuous relaxation of the problem. One is interested in finding a halfspace in the space of finite support functions $y: G \to \R$ that contains $I_\bve(G)$ but not $y_0$. Any such halfspace is described by an inequality of the form $\sum_{\x\in G}\pi(\x)y(\x) \geq 1$, where the function $\pi:G\to \R$ (which is not necessarily of finite support) gives the coefficients of the inequality. We use the notation $H_\pi$ to denote this halfspace, and call $\pi$ a {\em valid function} whenever $I_\bve(G)\subseteq H_\pi$.


Most of the literature has focused on the family of valid functions $\pi\ge 0$ and recent results \cite{basu2017structure} justify this nonnegativity restriction on $\pi$ for groups $G$ that are most relevant for integer programming. We will work with nonnegative valid functions in this paper. A valid function $\pi:G\to\R_+$ is {\em minimal} if every valid function $\tilde\pi:G\to\R_+$ such that $\tilde\pi\le\pi$ satisfies $\tilde\pi=\pi$.
We denote by $\mathcal M_\bve(G)$ the set of minimal functions.  The reason for restricting to minimal functions is that if $\pi:G\to\R_+$ is a valid function, then there always exists  a minimal $\pi'\in \mathcal M_\bve(G)$ such that $I_\bve(G)\subseteq H_{\pi'}\cap\R_+^G\subseteq H_\pi\cap\R_+^G$. For any group, minimal functions are characterized by the following theorem \cite{infinite2}.

\begin{theorem}\label{th:minimal}
A function $\pi:G\to\R_+$ is a minimal function if and only if it is subadditive (i.e., $\pi(\x)+\pi(\y)\ge\pi(\x+\y)$ for every $\x,\y\in G$), $\pi(\0) = 0$, and $\pi(\x) + \pi(\bve- \x) = 1$ for every $\x\in G$ (this property is known as the {\em symmetry property}).
\end{theorem}

A valid function $\pi:G\to\R_+$ is \emph{extreme} if $\pi_1=\pi_2$ for every pair of valid functions $\pi_1,\pi_2:G\to\R_+$ such that $\pi=\frac{1}{2}\pi_1+\frac{1}{2}\pi_2$. It is well known that the set of extreme functions is a subset of $\mathcal M_\bve(G)$. In fact, if $G$ is a finite group, the set $\mathcal M_\bve(G)$ is a polytope defined as\begin{equation}\label{eq:polytope}
\left\{\pi\in\R^G_+:\begin{array}{rcll}\pi(\0) & = &0&\\ \pi(\x)+\pi(\y) & \ge &\pi(\x+\y)&\forall\x,\y\in G \\ \pi(\x)+\pi(\bve- \x)& = &1&\forall\x\in G\end{array}\right\},
\end{equation}
and the extreme functions are precisely the vertices of this polytope.
\bigskip

\begin{example}\label{ex:minimal}
We give some well-known examples of minimal functions.
\begin{enumerate}
\item Let $G=\T^n$ and $\bve\in \T^n\setminus \{0\}$, and let $i \in \{1, \ldots, n\}$ be any index such that $\bve_i\ne0$. The {\em Gomory function $\GMI^n_{\bve_i}: \T^n \to \R_+$} is defined as $\GMI^n_{\bve_i}(\x) := \frac{\x_i}{\bve_i}$ if $\x_i \leq \bve_i$ and $\GMI^n_{\bve_i}(\x) := \frac{1-\x_i}{1-\bve_i}$ otherwise, where on the right hand sides the coordinates $\x_i, \bve_i$ are interpreted as real numbers in the interval $[0,1)$. It is well-known that $\GMI^n_{\bve_i}$ is an extreme function. This function is used to derive the well-known GMI cut in integer programming literature.

\item For $G = \Z/q\Z$, where $q$ is any natural number, and $b \in G\setminus \{0\}$, the function $\GOM^q_{b}:G\to \R_+$ defined as $\GOM^q_{b}(x) := \frac{x}{b}$ when $x \leq b$ and $\GOM^q_{b}(x) := \frac{q-x}{q-b}$ when $x > b$ is extreme.
\item Let $G = \Z/q\Z$ where $q$ is prime. Let $b_1, b_2 \in G\setminus \{0\}$ and let $\pi$ be a valid function for $I_{b_1}(G)$. Consider the automorphism $\phi: G\to G$ that sends $b_2$ to $b_1$ (which exists and is unique because $q$ is prime). Then $\pi\circ\phi$ is a valid function for $I_{b_2}(G)$.
    Furthermore $\pi$ is  minimal (extreme) if and only if  $\pi\circ\phi$ is  minimal (extreme).
\end{enumerate}
\end{example}


\section{Criteria to evaluate the strength of a valid function}\label{sec-volume}

\subsection{The distance criterion}

Consider the case when $G$ is a finite group. Given a valid function $\pi$, the halfspace $H_\pi:=\left\{y \in \R^G: \sum_{\x\in G}\pi(\x)y(\x)\ge1\right\}$ cuts off a portion of the nonnegative orthant $\R^G_+$ that includes the origin. A possible measure one can choose to evaluate to strength of $\pi$ is the distance of the origin from the hyperplane $\left\{y \in \R^G: \sum_{\x\in G}\pi(\x)y(\x)=1\right\}$: the larger this distance, the better $\pi$ is. It can be shown that this distance is in fact given by $\frac{1}{\lvert \pi\rvert_2}$, where $\lvert \pi\rvert_2 := \sqrt{\sum_{\x\in G}\pi(\x)^2}$ is the standard $\ell_2$ norm of the vector of coefficients given by $\pi$. Thus, one looks for $\pi$ that minimizes the $\ell_2$ norm. Generalizing, one could also look for a function $\pi$ that minimizes the $\ell_p$ norm for $p\geq 1$.

If $G$ is a group of infinite order, it is important to have a measure on the group $G$, with respect to which one can integrate so that expressions involving sums in the finite group case can be replaced by integrals. This is why we will consider compact, topological groups with the Haar measure. Since we would like to integrate functions defined in a topological group, we will have to restrict our attention to measureable functions (with respect to the Haar measure). Without further comment, we will assume that all functions $\pi$ that we consider are measurable; in particular, $\mathcal M_\bve(G)$ is to be thought of as the subset of {\em measurable} minimal functions from this point on. (For groups of finite order this is no restriction because all functions are measurable.) 
Under this assumption, for every $p\ge1$ the $\ell_p$ norm can be generalized to the standard $L_p$ norm $|\pi|_p := (\int |\pi|^p d\mu)^{1/p}$, where $\mu$ is the corresponding Haar measure.

The following result shows that the above optimization problem, i.e., maximizing the distance of the origin from the hyperplane defined by $\pi$, has a trivial solution: the $L_p$ norms are minimized by the function which takes value $1/2$ almost everywhere. The result holds for any group -- finite or infinite. Note that in the following statement we restrict to functions $\pi\in\mathcal M_\bve(G)$, as it is easy to verify that every optimal solution to this problem is a minimal function.

\begin{theorem}\label{thm:ell-p} Let $G$ be any (finite or infinite) group and $\bve \in G\setminus\{\0\}$. Let $p \geq 1$. Then $\inf\{|\pi|_p : \pi \in \mathcal M_\bve(G)\} = 1/2$ and is attained by the function which takes value $1/2$ everywhere except for $\0$ and $\bve$, where it takes value $0$ and $1$ respectively.
\end{theorem}

\proof{Proof.}
The proof of this theorem follows directly from the symmetry of $\pi$ (see Theorem \ref{th:minimal}). Indeed, since $\pi$ is symmetric, $D := \pi^{-1}([0,1/2))$ and $D' := \pi^{-1}((1/2,1])$ have equal measure. Moreover,
\[|\pi|^p_p = \int_D |\pi|^pd\mu + \int_{D'} |\pi|^pd\mu + \int_{\pi^{-1}(\{1/2\})}(1/2)^p d\mu.\]
Since $\mu(D)=\mu(D')$ and $\pi$ is symmetric, we obtain that
\[\begin{split}
|\pi|_p^p &= \int_D |\pi|^pd\mu + \int_{D} |1-\pi|^pd\mu + \int_{\pi^{-1}(\{1/2\})}(1/2)^p d\mu \\
&= \int_D (|\pi|^p + |1-\pi|^p)d\mu + \int_{\pi^{-1}(\{1/2\})}(1/2)^p d\mu.
\end{split}\]
For any real number $0 \leq a\leq 1$, we have that $a^p + (1-a)^p \geq 2(1/2)^p$ by convexity of the function $x\mapsto x^p$. Since $\pi$ is bounded between 0 and 1 by Theorem \ref{th:minimal}, we obtain $|\pi|_p \geq 1/2$.
\endproof
\medskip

The function which takes value $1/2$ everywhere (except $0$ and $\bve$) gives rise to cutting planes in integer programming that have been called MD2 cuts in the literature~\cite{rubin1972strengthened} and were obtained as a strengthening of the so-called {\em Dantzig cut}, which comes from the valid function which takes value $1$ everywhere~\cite{dantzig1959note}. In~\cite{rubin1972strengthened}, MD2 cuts are shown to yield a finite cutting plane algorithm for pure integer programming problems, which is based on earlier results from~\cite{bowman1970finiteness}.

It is interesting to observe that the above optimal function which takes value $1/2$ almost everywhere is not an extreme function.
Below we suggest another criterion to evaluate the strength of $\pi$, and we will see that it is optimized by an extreme function.

\subsection{The volume criterion}\label{sec:volume}

Consider first the case when $G$ is a finite cyclic group of prime order, i.e., $G=\Z/q\Z$ with $q$ prime.
Take $\pi \in \mathcal M_b(G)$. We can easily verify that $\pi(x)>0$ for every $x\in G\setminus\{0\}$. Indeed, for any $x\in G\setminus\{0\}$, there exists $k\in\N$ such that $kx=b$ because $q$ is prime. Since by Theorem \ref{th:minimal} $\pi$ is subadditive and $\pi(b)=1$, we have $k\pi(x)\geq \pi(kx)=\pi(b)=1$, hence $\pi(x)\geq 1/k$. 

Denote by $e^{x}$, $x\in G$, the unit vectors in $\R^G$. The above observation implies that the halfspace $H_{\pi}$ is parallel to $e^0$ and cuts off a simplex from the $(q-1)$-dimensional nonnegative orthant $\R^{G\setminus\{0\}}_+$. This simplex is given by $\conv\left(\{0\} \cup \left\{\frac{1}{\pi(x)}e^{x}\right\}_{x\in G\setminus\{0\}}\right)$. Consequently, its $(q-1)$-dimensional volume is given by $\frac1{(q-1)!}\prod_{x\in G\setminus\{0\}}\frac{1}{\pi(x)}$. We consider this volume as a measure of $\pi$: the higher this volume, the ``better" the halfspace given by $\pi$ is. Thus, one looks for $\pi \in \mathcal M_b(G)$ that minimizes $\prod_{x\in G\setminus\{0\}}\pi(x)$.
Since valid functions cut off the function $y$ defined by $y(x)=0$ for all $x\in G$, and $I_b(G)$ is contained in the nonnegative orthant, this provides a justification of the volume criterion as a measure of the strength of $\pi$.

We remark that, unlike the distance criterion, the function $\pi\mapsto\prod_{x\in G\setminus\{\0\}}\pi(x)$ is strictly log-concave, and therefore every minimizer $\pi$ is a vertex of the polytope in \eqref{eq:polytope}, i.e., $\pi$ is an extreme function.\medskip


The above definitions were made for a finite group $G$. We now show how to extend the volume measure defined above to infinite groups.

Let $G$ be a compact topological group with Haar measure $\mu$. Take $\pi \in \mathcal M_\bve(G)$. Let us initially assume that $\pi(\x)>0$ for every $\x\in G\setminus\{\0\}$. We denote by $\R^{(G)}$ the space of finite support functions $y:G\to\R$. Note that a basis of the vector space $\R^{(G)}$ is given by the family of functions $e^{\x}$, $\x\in G$, where $e^\x$ is the function which takes value 1 on $\x$ and 0 elsewhere. Similar to the finite cyclic group case, the halfspace $H_\pi$ is parallel to $e^\0$ and cuts off a convex set (let us call it a ``simplex'') from the set $\R^{(G)}_+\cap\{y:y(\0)=0\}$. This ``simplex'' is given by $\conv\left(\{0\} \cup \left\{\frac{1}{\pi(\x)}e^{\x}\right\}_{\x\in G\setminus\{\0\}}\right)$. However, we cannot compute the volume of this set, as it is an infinite dimensional object. To overcome this difficulty, we observe that in the finite cyclic group case maximizing the volume of the simplex is equivalent to maximizing its average side length (geometric mean). Therefore, in the infinite group case, we look at the geometric mean of the sides of the ``simplex'', i.e., the geometric mean of the function $1/\pi$, which is defined as
\[\exp\left(\frac1{\mu(G\setminus\{\0\})}\int_{G\setminus\{\0\}}\ln(1/\pi)d\mu\right)=\exp\left(\int_G\ln(1/\pi)d\mu\right).\]
(The equality holds because the Haar measure satisfies the properties $\mu(G)=1$ and $\mu(\{\x\})=\mu(\{\y\})$ for every $\x,\y\in G$, thus $\mu(\{\x\})=0$ for every $\x\in G$ because $G$ is infinite; in particular, $\mu(\{\0\})=0$ and $\mu(G\setminus\{\0\})=1$.)
Equivalently, we will minimize $\int_G\ln(\pi)d\mu$.

While we have motivated this formula for functions that are strictly positive everywhere, for minimal functions, this restriction is not necessary for the integral to make sense. Indeed, we will be concerned with integrals of functions of the form $\ln(\pi)$, where $\pi$ is bounded between $0$ and $1$. This means that $-\ln(\pi)$ is a nonnegative function taking values in the extended reals, i.e., it could take the value $+\infty$ at some points where $\pi$ equals $0$. For nonnegative, extended real valued functions, integrals are always defined but may equal $+\infty$. Below, we will say a nonnegative extended real valued function is {\em integrable} if the integral is finite.

\subsection{Our results}

With the above setup, we state our main results. The first result is about the volume measure for finite cyclic groups of prime order $q$, and shows that the maximum volume is cut off by appropriate automorphisms of the functions $\GOM^q_{q-1}$ and $\GOM^q_1$ (see part 3. of Example~\ref{ex:minimal} for these functions).

\begin{theorem}\label{thm:finite-volume} Let $G = \Z/q\Z$ where $q$ is prime, and let $b \in G\setminus \{0\}$. Then
\[\inf\left\{\prod_{x\in G\setminus\{0\}}\pi(x): \pi \in \mathcal M_b(G)\right\} = \frac{(q-1)!}{(q-1)^{q-1}}.\]
Moreover, this infimum is attained uniquely by the function $\GOM^q_{q-1}\circ \phi$, where $\phi: G\to G$ is the automorphism that sends $b$ to $q-1$.
\end{theorem}

Note that $\phi$ exists and is unique because $q$ is prime, and $\phi$ is equal to the multiplication by $-b^{-1}$, where $b^{-1}$ is the inverse of $b$ in the field $\Z/q\Z$. Note also that $\GOM^q_{q-1}\circ \phi$ = $\GOM^q_{1}\circ \phi'$ where $\phi':G \to G$ is the automorphism that sends $b$ to $1$.

For infinite, connected groups, we show that the average side length of the ``simplex'' cut off by an optimal inequality is $e$ (the base of the natural logarithm). Moreover, in the usual case $G=\T^n$, the Gomory function $\GMI^n_{\bve_i}$ defined in Example \ref{ex:minimal} is an optimal inequality.

\begin{theorem}\label{thm:infinite-volume} Let $G$ be any compact, connected topological group with Haar measure $\mu$, and $\bve \in G\setminus\{\0\}$. Then for every $\pi \in \mathcal M_\bve(G)$, $-\ln(\pi)$ is integrable and
\[\inf\left\{\int_G \ln (\pi)d\mu : \pi \in \mathcal M_\bve(G)\right\} = -1.\]
If $G = \T^n$, the infimum is attained by $\GMI^n_{\bve_i}$ for any index $i$ such that $\bve_i\ne0$.
\end{theorem}

We remark that although the two theorems above may seem of different nature, their meaning is essentially the same. Indeed, if in the discrete case we minimize the natural logarithm of the geometric mean of the side lengths of the simplex instead of its volume, Theorem \ref{thm:finite-volume} yields
\[\begin{split}
\inf\left\{\ln\bigg(\bigg(\prod_{x\in G\setminus\{0\}}\pi(x)\bigg)^{\frac{1}{q-1}}\bigg): \pi \in \mathcal M_b(G)\right\}=& \inf\left\{\frac1{q-1}\ln\bigg(\prod_{x\in G\setminus\{0\}}\pi(x)\bigg): \pi \in \mathcal M_b(G)\right\} \\
& = \frac1{q-1}\ln\left(\frac{(q-1)!}{(q-1)^{q-1}}\right)\\
&= \frac1{q-1}\left(\ln((q-1)!)-\ln((q-1)^{q-1})\right)\\
& =\frac1{q-1}\left(-(q-1)+o(q)\right),
\end{split}
\]
where the last equality follows from the fact that $\ln((q-1)^{q-1})=(q-1)\ln(q-1)$ and, by Stirling's approximation, $\ln((q-1)!)=(q-1)\ln(q-1)-(q-1)+o(q)$.
If we take the limit for $q\to\infty$, the above infimum tends to $-1$. This shows that Theorem \ref{thm:finite-volume} has a deep similarity with Theorem \ref{thm:infinite-volume}. In fact, a part of the proof of Theorem~\ref{thm:infinite-volume} involves using Theorem \ref{thm:finite-volume} and applying a limit argument on approximating Riemann sums for a Riemann integral.

The proofs of Theorems~\ref{thm:finite-volume} and~\ref{thm:infinite-volume} form the content of the rest of the paper.

\section{Finite cyclic groups of prime order: Proof of Theorem~\ref{thm:finite-volume}}

The main tool behind the proof of Theorem~\ref{thm:finite-volume} is a rearrangement idea of the function values which preserves the properties of subadditivity and symmetry, and makes the ``rearranged'' function nondecreasing. The contents of this idea are summarized in Theorem~\ref{thm:rearrange-finite}. One then shows that within the family of nondecreasing, subadditive and symmetric functions, any extreme function that minimizes the volume measure is the Gomory function; this is shown in Lemma~\ref{lem:extreme-is-Gomory}. Putting these two results together gives us Theorem~\ref{thm:finite-volume} (see the last proof of this section).

In the following, when we say that a function $\pi:\Z/q\Z\to\R_+$ is non decreasing, we mean that $\pi(x)\le \pi(y)$ for every $x,y\in\Z/q\Z$ with $x\le y$, where $x$ and $y$ are viewed as numbers in $\{0,\dots,q-1\}$ with the standard ordering.

\begin{theorem}\label{thm:rearrange-finite}
Let $G = \Z/q\Z$ with $q$ prime, and let $\pi: G\to \R_+$ be a subadditive function on $G$ such that $\pi(0) =0$ and $\pi$ is not identically $0$. Then $\hat \pi: G \to \R_+$ defined as
\begin{equation}\label{eq:rearrange-2}
\hat \pi(x) := \min \{\alpha \geq 0: |\pi^{-1}((0,\alpha])| \geq x\}\quad \forall x \in G
\end{equation}
is finite-valued, subadditive on $G$ and nondecreasing. (On the right-hand side, $x$ is viewed as a number in $\{0,\dots,q-1\}$.) Moreover, $|\pi^{-1}(\{\beta\})| = |\hat \pi^{-1}(\{\beta\})|$ for all $\beta > 0$. Finally, if there exists some $b \in G$ such that $\pi(x) + \pi(b- x) = 1$ for all $x\in G$, then $\hat \pi(x) + \hat \pi(q-1-x) = 1$ for all $x\in G$.
\end{theorem}

\proof{Proof.}
The argument given at the beginning of Subsection \ref{sec:volume} shows that $\pi(x) > 0$ for all $x\neq 0$, since $\pi$ is subadditive and not identically $0$, and since $q$ is prime. Thus, for every $x \in \{0, \ldots, q-1\}$, there exists $\alpha \geq 0$ such that $|\pi^{-1}((0,\alpha])| \geq x$. So the minimum in~\eqref{eq:rearrange-2} is taken over a nonempty set. Since $G$ is finite, this set has a minimum, which is attained at some  $\alpha\in\pi(G)$. This shows that the minimum is well-defined.
The nondecreasing property follows from the definition.

For the proof of subadditivity we need to use the Cauchy--Davenport theorem~\cite{nathanson}, which states that for any nonempty subsets $A, B\subseteq G$, we have the inequality
$|A+B|\geq \min\{q,\, |A|+|B|-1\}$  (where $A+B:=\{\a+\bve:a\in A,\,b\in B\}$). Consequently, any nonempty $A, B\subseteq G$ such that $0 \not\in B$ satisfy $|(A+B)\cup A|\geq \min\{q,\, |A|+|B|\}$: apply the Cauchy--Davenport theorem to $A, B\cup\{0\}$.

Consider $x_1, x_2\in G$ and let $\alpha_i = \hat \pi(x_i)$, $i=1,2$. Therefore, $|\pi^{-1}((0, \alpha_i])| \geq x_i$, $i=1,2$. Moreover, by subadditivity of $\pi$, $$\left(\pi^{-1}((0,\alpha_1]) + \pi^{-1}((0,\alpha_2])\right) \cup \pi^{-1}((0,\alpha_1])\subseteq \pi^{-1}((0,\alpha_1+\alpha_2]).$$ Therefore, we obtain
$$
\begin{array}{rcl}
|\pi^{-1}((0,\alpha_1 + \alpha_2])| & \geq & \left|\left(\pi^{-1}((0,\alpha_1]) + \pi^{-1}((0,\alpha_2])\right) \cup \pi^{-1}((0,\alpha_1])\right| \\
&\geq &\min\{q,\, |\pi^{-1}((0,\alpha_1])| + |\pi^{-1}((0,\alpha_2])|\} \\
& \geq & x_1 + x_2
\end{array}
$$
where we use the fact that $\pi(0) =0$, and so $0 \not\in \pi^{-1}((0,\alpha_2]))$, and also used the Cauchy-Davenport theorem for the second inequality. Note that the last term $x_1+x_2$ is a sum in $G$ and is viewed as a real number, thus $x_1+x_2<q$. This establishes the subadditivity of $\hat \pi$.
%
%

We next confirm that $|\pi^{-1}(\{\beta\})| = |\hat \pi^{-1}(\{\beta\})|$ for all $\beta > 0$. To show this, it suffices to prove that $|\pi^{-1}((0,\beta])| = |\hat \pi^{-1}((0,\beta])|$ for all $\beta > 0$ . Observe that for any $t \in G$,
$$\hat \pi(t) \leq \beta \;\;\Leftrightarrow\;\; \min\{\alpha \geq 0: |\pi^{-1}((0,\alpha])| \geq t\} \leq \beta\;\; \Leftrightarrow\;\; |\pi^{-1}((0,\beta])| \geq t.$$
Therefore, since $\hat \pi$ is nondecreasing and $\hat\pi(0)=0$,
$$|\hat \pi^{-1}((0,\beta])| = \max\{t \in G: \hat \pi(t) \leq \beta\} = \max\{t \in G: |\pi^{-1}((0,\beta])| \geq t\} = |\pi^{-1}((0,\beta])|.$$

We finally verify that the symmetry condition is preserved. Consider any $x \in G$ and let $\alpha = \hat \pi(x)$. By definition, $|\pi^{-1}((0,\alpha])| \geq x$. By symmetry of $\pi$, $t \in \pi^{-1}((0,\alpha))$ if and only if $b- t \in \pi^{-1}((1-\alpha,1))$. 
We now observe that $\alpha = \hat \pi(x)$ implies that
$$\begin{array}{rrcl}& x& \geq & |\pi^{-1}((0,\alpha))| \\
\Rightarrow& x & \geq & |\pi^{-1}((1-\alpha, 1))| \\
\Rightarrow& x & \geq & q-1-|\pi^{-1}((0,1-\alpha])| \\
\Rightarrow& |\pi^{-1}((0,1-\alpha])| & \geq &q-1-x.
\end{array}$$
Thus, by definition of $\hat \pi$, $\hat \pi(q-1-x) \leq 1-\alpha$. Therefore, we have $\hat \pi(x) + \hat \pi(q-1-x) \leq 1$. By subadditivity of $\hat \pi$, $\hat \pi(x) + \hat \pi(q-1-x) \geq \hat \pi(q-1)$. We finally observe that $\hat \pi(q-1) = 1$ because $\pi(b) = 1$ by symmetry of $\pi$ and $\pi(x) \leq 1$ for all $x\in G$. Therefore, $\hat \pi(x) + \hat \pi(q-1-x) = 1$.\endproof
\medskip

\begin{lemma}\label{lem:rearrange-GOM-q-1}
Let $G = \Z/q\Z$ where $q$ is prime, and let $b \in G\setminus \{0\}$. Then the function $h= \GOM^q_{q-1}\circ \phi$, where $\phi: G\to G$ is the automorphism that sends $b$ to $q-1$, satisfies $\hat h(x) = \frac{x}{q-1}$ for every $x\in G$, where on the right hand side $x$ is viewed as a real number.
\end{lemma}

\proof{Proof.} Since $\GOM^q_{q-1}$ takes all the distinct values in $\left\{0, \frac{1}{q-1}, \frac{2}{q-1}, \ldots, \frac{q-2}{q-1}, 1\right\}$, so does the function $h$. Since the rearrangement $\hat h$ is nondecreasing, we obtain the result.\endproof
\medskip

\begin{lemma}\label{lem:extreme-is-Gomory}
Let $G = \Z/q\Z$ where $q$ is prime. If $\pi\in\mathcal M_{q-1}(G)$ is an extreme function and is nondecreasing, then $\pi=\GOM^q_{q-1}$.
\end{lemma}

\proof{Proof.}
Let $x,y\in G$ and assume that $x+y\ge q$ (where $x,y,q$ are viewed as real numbers). Then $x+y$, taken modulo $q$, is smaller than $x$, and thus $\pi(x+y)\le\pi(x)$, as $\pi$ is nondecreasing. Furthermore, the condition $x+y\ge q$ implies $y\ne0$, and thus $\pi(y)>0$. This implies that $\pi(x)+\pi(y)>\pi(x+y)$.

By the above observation, the number $\gamma:=\min\{\pi(x)+\pi(y)-\pi(x+y):x,y\in G,\,x+y\ge q\}$ is strictly positive. Then there exists $\lambda\in(0,1)$ such that $\lambda\le\gamma(q-1)/q$ and $\lambda\le\pi(x)/x$ for every $x\in G\setminus\{0\}$ (where the denominator $x$ is interpreted as a real number).
For $x\in G$, define
\[\tilde\pi(x)=\frac{\pi(x)-\lambda g(x)}{1-\lambda},\]
where $g=\GOM^q_{q-1}$, i.e., $g(x)=\frac{x}{q-1}$ for every $x\in G$.
Note that $\pi=\lambda g+(1-\lambda)\tilde\pi$. Therefore, if we show that $\tilde\pi\in\mathcal M_{q-1}(G)$, then by the extremality of $\pi$ we have that $\pi=\tilde\pi=g$, and the proof is complete.

It remains to prove that $\tilde\pi\in\mathcal M_{q-1}(G)$. We do so by showing that $\tilde\pi$ fulfills the conditions of Theorem \ref{th:minimal}.
It is clear by definition that $\tilde\pi(0)=0$. Furthermore, since $\lambda\le\pi(x)/x$ for every $x\in G\setminus\{0\}$, we have $\tilde\pi(x)\ge0$ for every $x\in G$.
Symmetry holds because, by symmetry of $\pi$ and $g$, for $x\in G$ we have
\[\tilde\pi(x)+\tilde\pi(q-1-x)=\frac{\pi(x)+\pi(q-1-x)-\lambda(g(x)+g(q-1-x))}{1-\lambda}=\frac{1-\lambda}{1-\lambda}=1.\]
We finally check that $\tilde\pi$ is subadditive. If $x,y\in G$ satisfy $x+y<q$, then, by the subadditivity of $\pi$ and by definition of $g$,
\[\tilde\pi(x)+\tilde\pi(y)=\frac{\pi(x)+\pi(y)-\lambda(g(x)+g(y))}{1-\lambda}\ge\frac{\pi(x+y)-\lambda g(x+y)}{1-\lambda}=\tilde\pi(x+y).\]
If, on the contrary, $x+y\ge q$, then
\[g(x+y)=\frac{x+y-q}{q-1}=g(x)+g(y)-\frac{q}{q-1}\]
and thus
\[\begin{split}
\tilde\pi(x)+\tilde\pi(y) &=\frac{\pi(x)+\pi(y)-\lambda(g(x)+g(y))}{1-\lambda}\\
&= \frac{(\pi(x)+\pi(y)-\pi(x+y)-\lambda q/(q-1))+\pi(x+y)-\lambda g(x+y)}{1-\lambda}\\
&\ge \frac{\pi(x+y)-\lambda g(x+y)}{1-\lambda}=\tilde\pi(x+y),
\end{split}\]
where the inequality holds by definition of $\lambda$.\endproof
\medskip

\proof{Proof of Theorem~\ref{thm:finite-volume}.}
Let $G = \Z/q\Z$ where $q$ is a prime number, and let $b \in G\setminus \{0\}$.
We want to show that
\[\inf\left\{\prod_{x\in G\setminus\{0\}}\pi(x): \pi \in \mathcal M_b(G)\right\} = \frac{(q-1)!}{(q-1)^{q-1}}\]
and the infimum is attained uniquely for $\pi=\GOM^q_{q-1}\circ \phi$, where $\phi: G\to G$ is the automorphism that sends $b$ to $q-1$.
We denote by $\mathcal P_b(G)$ the above optimization problem.

Recall the discussion after Theorem \ref{th:minimal}, where $\mathcal M_b(G)$ is viewed as a polytope in $\R^G$ defined by the linear inequalities in~\eqref{eq:polytope}. This implies that the optimization problem $\mathcal P_b(G)$ has an optimal solution, as the objective function is continuous.
Furthermore, since the objective function is strictly log-concave, every minimizer is an extreme point of the polytope $\mathcal M_b(G)$. In other words, every minimizer is an extreme function for $\mathcal M_b(G)$.

Assume first that $b=q-1$ and let $\pi$ be an optimal solution to $\mathcal P_{q-1}(G)$. By Theorem~\ref{thm:rearrange-finite}, the function $\hat\pi$ defined in \eqref{eq:rearrange-2} belongs to $\mathcal M_{q-1}(G)$ and has the same objective value as $\pi$, thus it is also an optimal solution to $\mathcal P_{q-1}(G)$. This implies that $\hat\pi$ is an extreme function for $\mathcal M_{q-1}(G)$ and therefore, by Lemma \ref{lem:extreme-is-Gomory}, $\hat\pi=\GOM_{q-1}^q$ and $\prod_{x\in G\setminus\{0\}}\hat\pi(x)=\frac{(q-1)!}{(q-1)^{q-1}}$.

Assume now $b\ne q-1$.
By Lemma~\ref{lem:rearrange-GOM-q-1}, the function $h=\GOM^q_{q-1}\circ \phi$, where $\phi: G\to G$ is the automorphism that sends $b$ to $q-1$, satisfies $\prod_{x\in G\setminus\{0\}}h(x)=\frac{(q-1)!}{(q-1)^{q-1}}$. Thus the minimum value of $\mathcal P_b(G)$ is at most $\frac{(q-1)!}{(q-1)^{q-1}}$.
Let $\pi$ be an optimal solution to $\mathcal P_b(G)$ and assume that the optimal value is smaller than $\frac{(q-1)!}{(q-1)^{q-1}}$. By Theorem~\ref{thm:rearrange-finite}, the function $\hat\pi$ defined in \eqref{eq:rearrange-2} belongs to $\mathcal M_{q-1}(G)$ and has the same objective value as $\pi$, contradicting the fact that the optimal value of $\mathcal P_{q-1}(G)$ is $\frac{(q-1)!}{(q-1)^{q-1}}$.

It remains to prove that, for any $b\in G\setminus\{0\}$, $\GOM_{q-1}^q\circ\phi$ is the unique optimal solution to $\mathcal P_b(G)$, where $\phi: G\to G$ is the automorphism that sends $b$ to $q-1$. The above arguments show that if $\pi$ is an optimal solution to $\mathcal P_b(G)$, then $\hat\pi=\GOM_{q-1}^q\circ\phi$, and thus each of $\pi$ and $\hat\pi$ takes all the distinct values in $\left\{0,\frac{1}{q-1},\frac{2}{q-1},\dots,1\right\}$. To conclude the proof, we show that any function $\pi\in\mathcal M_b(G)$ that takes all the distinct values in $\left\{0,\frac{1}{q-1},\frac{2}{q-1},\dots,1\right\}$ coincides with $\GOM_{q-1}^q\circ\phi$.

Let $x\in G\setminus\{0\}$ be such that $\pi(x)=\frac1{q-1}$. By subadditivity, $\pi(2x)\le2\pi(x)=\frac2{q-1}$, and thus $\pi(2x)=\frac2{q-1}$. Iterating this argument yields $\pi(kx)=\frac{k}{q-1}$ for every $k\in\{0,\dots,q-1\}$.
This means that $\pi=\GOM_{q-1}^q\circ\psi$, where $\psi: G\to G$ is the automorphism that sends 1 to $x$.
Now let $\bar k$ be the unique number in $\{0,\dots,q-1\}$ such that $\bar kx=b$ in $G$. (Such a number exists because $q$ is prime.) Then $1=\pi(b)=\pi(\bar kx)=\frac{\bar k}{q-1}$, thus $\bar k=q-1$. Then $\psi$ sends $b$ to $q-1$ and therefore $\psi=\phi$.\endproof
\medskip

\section{Infinite, connected groups: Proof of Theorem~\ref{thm:infinite-volume}}

In order to prove Theorem~\ref{thm:infinite-volume}, we will show an analogue of Theorem~\ref{thm:rearrange-finite}. The reader will notice that Theorem~\ref{thm:rearrange-finite} permutes the finitely many function values into non-decreasing order, i.e., sorts the function values. We will try to do the same for infinite, connected groups. Of course, now we have to do some sort of ``continuous" sorting; we will now preserve integrals as opposed to exact function values. More precisely, we use a rearrangement idea which preserves the properties of subadditivity, symmetry and the values of integrals, and makes the ``rearranged'' function nondecreasing. We then prove Theorem \ref{thm:infinite-volume} by applying a limit argument to Theorem~\ref{thm:finite-volume}.

We now begin with the appropriate definitions and results needed to make the rearrangement idea concrete. In the following, given $x\in\T^1$, when we say ``$x$ viewed as a real number'' we refer to the unique representative of $x$ in the interval $[0,1)$.

\begin{definition} Let $G$ be any compact, topological group with Haar measure $\mu$ on it. For any nonnegative function $\pi: G\to \R_+$, define $\hat \pi: \T^1 \to \R_+$ as follows:
\begin{equation}\label{eq:rearrange}
\hat \pi(x) := \inf\{\alpha \geq 0: \mu(\pi^{-1}([0,\alpha])) \geq x\}\quad \forall x \in \T^1,
\end{equation}
where the right hand side of the inequality $\mu(\pi^{-1}([0,\alpha])) \geq x$ is viewed as a real number.
\end{definition}


\begin{remark} We note that for any $x \in [0,1)$, the set $\{\alpha \geq 0: \mu(\pi^{-1}([0,\alpha])) \geq x\}$ is nonempty, and so the infimum in~\eqref{eq:rearrange} is a well-defined real number. Indeed, since $\pi$ is nonnegative, we have $\bigcup_{n \in \N}\pi^{-1}([0,n]) = \pi^{-1}([0,\infty)) = G$. Then, by continuity of measure, $\lim_{n\to \infty}\mu(\pi^{-1}([0,n])) = \mu(G) = 1$. Thus, for any $x \in [0,1)$, there must exist some $n \in \N$ such that $\mu(\pi^{-1}([0,n])) \geq x$.
\end{remark}

\begin{lemma}\label{lem:two-prop} Let $G$ be any compact, connected topological group with Haar measure $\mu$ on it, and let $\pi: G\to \R_+$. The following hold.
\begin{enumerate}
\item\label{two-prop:1} The function $\alpha \mapsto \mu(\pi^{-1}([0,\alpha]))$ is nondecreasing and right continuous.
\item\label{two-prop:2} $\hat \pi$ is nondecreasing, i.e., $x \leq y$ implies $\hat \pi(x) \leq \hat \pi(y)$, where $x,y\in \T^1$ are viewed as real numbers with the standard ordering.
\item\label{two-prop:3} Let $\bar x \in \T^1$ and let $\bar \alpha = \hat \pi(\bar x)$. Then $\mu(\pi^{-1}([0,\bar\alpha)))\leq \bar x \leq \mu(\pi^{-1}([0,\bar\alpha]))$, where $\bar x$ is viewed as a real number.
\item\label{two-prop:4} $\hat \pi$ is left continuous, i.e., for all $x \in \T^1\setminus\{0\}$, $\lim_{\epsilon \to 0^+} \hat \pi(x -  \epsilon) = \hat \pi(x)$.
\end{enumerate}
\end{lemma}

\proof{Proof.} 
For the first property, we observe that $[0,\alpha] = \bigcap_{\epsilon > 0} [0, \alpha + \epsilon]$, and therefore $\pi^{-1}([0,\alpha]) = \bigcap_{\epsilon > 0} \pi^{-1}([0, \alpha + \epsilon])$. By continuity of measure, $\lim_{\epsilon \to 0^+} \mu(\pi^{-1}([0, \alpha + \epsilon])) = \mu(\pi^{-1}([0,\alpha]))$, establishing property~\ref{two-prop:1}. The fact that the function is nondecreasing is clear from its definition.


The second property is clear from the definition of $\hat \pi$.

For the third property, we first show $\bar x\leq \mu(\pi^{-1}([0,\bar\alpha]))$. By definition of $\bar \alpha$, for every  $\epsilon>0$, $\mu(\pi^{-1}([0,\bar\alpha+\epsilon]))\geq \bar x$, therefore by part 1 we have $\mu(\pi^{-1}([0,\bar\alpha]))=\lim_{\epsilon\to 0^+}\mu(\pi^{-1}([0,\bar\alpha+\epsilon]))\geq \bar x$.
We show $\mu(\pi^{-1}([0,\bar\alpha))) \leq \bar x$. We have $[0,\bar\alpha)=\bigcup_{0<\epsilon<\bar\alpha}[0,\bar\alpha-\epsilon]$, hence $\pi^{-1}([0,\bar\alpha))=\bigcup_{0<\epsilon<\bar\alpha}\pi^{-1}([0,\bar\alpha-\epsilon])$. By definition of $\bar \alpha$,
$\mu(\pi^{-1}([0, \bar\alpha - \epsilon]))<\bar x$ for every $\epsilon\in(0,\bar\alpha)$. By continuity of measure, $\mu(\pi^{-1}([0,\bar \alpha)))=\lim_{\epsilon \to 0^+} \mu(\pi^{-1}([0, \bar \alpha - \epsilon])) \leq \bar x $.

For the fourth property, consider $x \in \T^1\setminus\{0\}$. For every $0 \leq \epsilon < x$, define $t_\epsilon := \hat \pi( x -  \epsilon)$. Since $\hat \pi$ is nondecreasing (property~\ref{two-prop:2} above), the function $\epsilon\mapsto t_\epsilon$ is nonincreasing and $t_\epsilon \leq t_0$ for all $0 < \epsilon <  x$. Therefore (see, e.g., \cite[Theorem 4.29]{babyrudin}), the limit $\bar t:= \lim_{\epsilon \to 0^+}t_\epsilon$ exists and $\bar t\leq t_0$. We need to show that $\bar t \geq t_0$. It suffices to prove that $\mu(\pi^{-1}([0,\bar t])) \geq x$, since $t_0 = \hat \pi( x) = \inf\{\alpha \geq 0:\mu(\pi^{-1}([0,\alpha])) \geq  x \}$. We observe that $[0,\bar t] \supseteq \bigcup_{0 < \epsilon < x}[0, t_\epsilon]$ 
since $\bar t\geq t_\epsilon$ for $\epsilon\in(0,x)$. By continuity of measure, $\mu(\pi^{-1}([0,\bar t])) \geq \lim_{\epsilon \to 0^+} \mu(\pi^{-1}([0,t_\epsilon])) \geq \lim_{\epsilon \to 0^+} ( x -  \epsilon) =  x$, where the second inequality follows from property 3.\endproof
\medskip

The next result makes use of Kemperman's theorem~\cite[Theorem 1.1]{kemperman1964products}, which, in our context, states that if $A$ and $B$ are nonempty subsets of a compact, connected group $G$, then $\mu(A+B)\ge\min\{1,\mu(A)+\mu(B)\}$. Kemperman's theorem can be seen as an analogue of the Cauchy--Davenport inequality used in the previous section. However, Kemperman's theorem only applies to connected groups, while the Cauchy--Davenport inequality only applies to finite cyclic groups, which are not connected. We highlight that the need to use Kemperman's theorem is the only reason for the connectedness assumption in Theorem \ref{thm:infinite-volume}.

\begin{theorem}\label{thm:rearrange} Let $G$ be any compact, connected topological group with Haar measure $\mu$ on it, and let $\pi: G\to \R_+$ be subadditive on $G$, with $\pi(0)=0$. Then $\hat \pi$ is subadditive on $\T^1$.\end{theorem}

\proof{Proof.} 
Consider $x_1, x_2 \in \T^1$. Let $\hat \pi(x_i) = \alpha_i$, $i=1,2$. It suffices to show that $\mu(\pi^{-1}([0,\alpha_1 + \alpha_2])) \geq x_1+x_2$, where the right hand side is a sum in $\T^1$ and its value is viewed as real number.

Subadditivity and nonnegativity of $\pi$ imply that $\pi^{-1}([0,\alpha_1 + \alpha_2]) \supseteq \pi^{-1}([0,\alpha_1]) + \pi^{-1}([0,\alpha_2])$. This implies:
\[\begin{split}
\mu\left(\pi^{-1}([0,\alpha_1 + \alpha_2])\right) & \ge \mu\left(\pi^{-1}([0,\alpha_1]) + \pi^{-1}([0,\alpha_2])\right)\\
& \geq \min\{1,\mu(\pi^{-1}([0,\alpha_1])) + \mu(\pi^{-1}([0,\alpha_2]))\} \\
& \geq x_1 + x_2,
\end{split}\]
where the second inequality follows from Kemperman's theorem (which can be applied because $\pi^{-1}([0,\alpha_i])\ne\emptyset$ for $i=1,2$, as $\pi(0)=0$), and the last inequality follows from the fact that $\mu(\pi^{-1}([0,\alpha_i])) \geq x_i$ because $\hat \pi(x_i) = \alpha_i$, for $i=1,2$ (property~\ref{two-prop:3} in Lemma~\ref{lem:two-prop}).\endproof
\medskip

\begin{lemma}\label{lem:rearrange-2} Let $G$ be any compact, connected topological group with Haar measure $\mu$ on it. Let $\bve \in G\setminus\{0\}$. Consider any $\pi:G\to\R_+$ such that $ \pi(\x) + \pi(\bve - \x) = 1$ for every $\x\in G$.
Then for any $x \in \T^1\setminus\{0\}$, $\hat \pi(x) + \hat \pi(-x) \leq 1$.
\end{lemma}

\proof{Proof.} Let $\hat \pi(x) = \alpha$.
Since $-x$, viewed as a real number, is $1-x$, we need to show that $\hat\pi(1-x)\le1-\alpha$.
For this, it suffices to prove that $\mu(\pi^{-1}([0,1-\alpha])) \geq 1-x$.
Define $S := \pi^{-1}([0,\alpha))$ and $S' := \{\bve- \x: \x \in S\} = \bve - S$. Since $\mu$ is the Haar measure on $G$, $\mu(S)=\mu(S')$. Moreover, as $\pi$ satisfies $\pi(\x) + \pi(\bve- \x) = 1$ for every $\x\in G$, we must have that $S' = \T^1\setminus \pi^{-1}([0,1-\alpha])$. Since $\hat \pi(x) = \alpha$, by property~\ref{two-prop:3} of Lemma~\ref{lem:two-prop} we obtain $x \geq \mu(S) = \mu(S') = 1 - \mu(\pi^{-1}([0,1-\alpha])).$
\endproof
\medskip

\begin{lemma}\label{lem:equi-meas} Let $G$ be any compact, connected topological group with Haar measure $\mu$ on it. For any nonnegative function $\pi: G\to \R_+$ and any $\beta \geq 0$, we have $$\ell(\hat \pi^{-1}([0,\beta])) = \mu(\pi^{-1}([0,\beta])),$$ where $\ell$ denotes the Haar measure on $\T^1$.
\end{lemma}

\proof{Proof.} By property \ref{two-prop:3} of Lemma~\ref{lem:two-prop}, for every $t \in \T^1$ and $\beta\in\R_+$, $\hat \pi(t) \leq \beta \;\; \Leftrightarrow\;\; \mu(\pi^{-1}([0,\beta])) \geq t$.
Therefore, since $\hat \pi$ is nondecreasing,
$$\ell(\hat \pi^{-1}([0,\beta])) = \sup\{t \geq 0: \hat \pi(t) \leq \beta\} = \sup\{t \geq 0: \mu(\pi^{-1}([0,\beta])) \geq t\} = \mu(\pi^{-1}([0,\beta])),$$ where the first supremum is taken over $t\in \T^1$ with the standard order on $\T^1$, and the second one is taken over $t \in \T^1$ viewed as a real number. Note that the first supremum is a well-defined real number because $\hat \pi(0)=0$.\endproof
\medskip

\begin{prop}\label{prop:preserve-int}  Let $G$ be any compact, connected topological group with Haar measure $\mu$ on it, and let $\pi: G\to [0,1]$ be subadditive on $G$. Let $\hat \pi$ be as defined in~\eqref{eq:rearrange}. Then $\int_G \ln(\pi)d\mu$ is finite if and only if $\int_{\T^1} \ln(\hat \pi)d\ell$ is finite, and in this case $$\int_G \ln(\pi)d\mu = \int_{\T^1} \ln(\hat \pi)d\ell,$$ where $\ell$ denotes the Haar measure on $\T^1$.
\end{prop}
\proof{Proof.} We use the so-called ``layer-cake representation'' of a nonnegative function (see, e.g.,~\cite[Chapter 1]{lieb2001analysis}), which states that for any nonnegative function $F$ defined on a measure space $(\Omega, \nu)$,
$$\int_\Omega Fd\nu = \int_0^\infty \nu(\{\omega \in \Omega: F(\omega) \geq t\})dt.$$
Therefore, since $\pi$ takes values in $[0,1]$, we have that $-\ln(\pi) \geq 0$ and we can write
\begin{equation}\label{eq:layer-1}
\begin{split}
\int_G -\ln(\pi)d\mu & = \int_0^\infty \mu(\{\x \in G: -\ln(\pi(\x)) \geq t\})dt \\
& = \int_0^\infty \mu(\{\x \in G: \pi(\x) \leq e^{-t}\})dt \\
& = \int_0^1 \frac{1}{s}\,\mu(\{\x \in G: \pi(\x) \leq s\})ds,
\end{split}
\end{equation}
where we use the change of variable $s = e^{-t}$. 
Since $\pi$ takes values in $[0,1]$, so does $\hat \pi$ and similarly we get \begin{equation}\label{eq:layer-2}\int_{\T^1} -\ln(\hat \pi)d\ell = \int_0^1 \frac{1}{s}\,\ell(\{x \in \T^1: \hat \pi(x) \leq s\})ds,\end{equation} 
Applying Lemma~\ref{lem:equi-meas} to~\eqref{eq:layer-1} and~\eqref{eq:layer-2} gives the desired result.\endproof
\medskip

We now consider another function derived from $\hat \pi$. Let $G$ be any compact, connected topological group with Haar measure $\mu$ on it. For any nonnegative function $\pi: G\to \R_+$, let $\hat \pi: \T^1 \to \R_+$ be defined as in \eqref{eq:rearrange}. Define

\begin{equation}\label{eq:right-limit-fn}
\bar \pi(x) := \lim_{\epsilon\to0^+} \hat \pi(x+\epsilon),
\end{equation}
which is well-defined because $\hat \pi$ is bounded from below and nondecreasing by property~\ref{two-prop:2} in Lemma~\ref{lem:two-prop}, and so the limit in the definition of $\bar \pi$ exists and is finite (see, e.g., \cite[Theorem 4.29]{babyrudin}).

\begin{lemma}\label{lem:right-limit-fn} Let $G$ be any compact, connected topological group with Haar measure $\mu$ on it and let $\bve \in G\setminus \{\0\}$. For any $\pi\in \mathcal M_\bve(G)$, let $\hat \pi, \bar \pi$ be defined as in \eqref{eq:rearrange} and~\eqref{eq:right-limit-fn}. Then $\bar \pi$ is nondecreasing and subadditive. Moreover, for any $x \in \T^1 \setminus \{0\}$, $\hat \pi(x) + \bar \pi(-x) = 1$.
\end{lemma}

\proof{Proof.}
We first note that since $\pi\in\mathcal M_\bve(G)$, $\pi$ satisfies the properties listed in Theorem \ref{th:minimal}.
Now, $\bar \pi$ is nondecreasing because $\hat \pi$ is nondecreasing by property~\ref{two-prop:2} in Lemma~\ref{lem:two-prop} . We check that $\bar \pi$ is subadditive. Consider $a,b \in \T^1$ and let $x = a + b$. Then
\[\begin{split}
\bar \pi(x) & = \lim_{\epsilon \to 0^+} \hat \pi(x+\epsilon) \\
& \leq  \lim_{\epsilon \to 0^+} (\hat \pi(a+\epsilon/2) + \hat \pi(b+\epsilon/2)) \qquad\textrm{by subadditivity of $\hat \pi$ (Theorem \ref{thm:rearrange})} \\
& =  \lim_{\epsilon \to 0^+} \hat \pi(a+\epsilon/2) + \lim_{\epsilon \to 0^+}\hat \pi(b+\epsilon/2) \\
& =  \bar \pi(a) + \bar \pi(b).
\end{split}\]
We now check that for any $x \in \T^1 \setminus \{0\}$, $\hat \pi(x) + \bar \pi(-x) = 1$. For every $0 < \epsilon < x$ we have that $\hat \pi(x -  \epsilon) + \hat \pi(-x + \epsilon) \leq 1$ by Lemma~\ref{lem:rearrange-2}. Taking the limit $\epsilon \to 0^+$, and using the fact that $\hat \pi$ is left continuous by property~\ref{two-prop:4} in Lemma~\ref{lem:two-prop}, we obtain that $\hat \pi(x) + \bar \pi(-x) \leq 1$. We show that the reverse inequality holds after establishing the following claim.
\begin{claim}\label{claim:alpha} For all $z \in \T^1$, $\bar \pi(z) =  \sup\{ t: \mu(\pi^{-1}([0,t))) \leq z\},$ where $z$ inside the supremum is viewed as a real number.
\end{claim}
\begin{cpf} Let $\bar t=\sup\{ t: \mu(\pi^{-1}([0,t))) \leq z\}$. We need to prove $\bar \pi(z)=\bar t$.
We first show $\bar \pi(z)\geq \bar t$. For every $\epsilon \in(0,1-z)$ and every $\alpha\in [0,\bar t)$, by definition of $\bar t$ we have $\mu(\pi^{-1}([0,\alpha]))<z+\epsilon$, which implies $\hat\pi(z+\epsilon)\geq\bar t$, and therefore $\bar \pi(z)=\lim_{\epsilon\to 0^+}\hat\pi(z+\epsilon)\geq \bar t$.

We  show $\bar \pi(z)\leq \bar t$. By property~\ref{two-prop:3} in Lemma~\ref{lem:two-prop}, $\mu(\pi^{-1}([0,\hat\pi(z))))\leq z$, which implies $\hat\pi(z)\leq \bar t$ by definition of $\bar t$. Again by property~\ref{two-prop:3} in Lemma~\ref{lem:two-prop}, it then follows that $z\leq \mu(\pi^{-1}([0,\bar t]))$.
Suppose that $z<\mu(\pi^{-1}([0,\bar t]))$, and let $\bar\epsilon=\mu(\pi^{-1}([0,\bar t]))-z>0$. Then, for every $\epsilon\in(0,\bar\epsilon]$, $\hat \pi(z+\epsilon)\leq \bar t$. This implies that $\bar\pi(z)=\lim_{\epsilon\to 0^+}\hat\pi(z+\epsilon)\leq \bar t$. Assume then that $\mu(\pi^{-1}([0,\bar t]))=z$. For every $\delta\in (0,1-\bar t)$, let $\epsilon(\delta)=\mu(\pi^{-1}([0,\bar t+\delta]))-z$. By definition of $\bar t$, $\epsilon(\delta)>0$ for every $\delta\in (0,1-\bar t)$. Furthermore,
$\lim_{\delta\to 0^+} \epsilon(\delta)=\lim_{\delta\to 0^+} (\mu(\pi^{-1}([0,\bar t+\delta]))-z)=\mu(\pi^{-1}([0,\bar t]))-z=0$, where the second equation follows from property~\ref{two-prop:1} of Proposition~\ref{lem:two-prop}. By definition, $\hat\pi(z+\epsilon(\delta))\leq \bar t+\delta$. Since $\epsilon(\delta)$ is strictly positive, nonincreasing, and converges to $0$ as $\delta\to 0^+$, it follows that
$\bar\pi(z)=\lim_{\epsilon\to 0^+}\hat \pi(z+\epsilon)=\lim_{\delta\to 0^+}\hat \pi(z+\epsilon(\delta))\leq \lim_{\delta\to 0^+} (\bar t+\delta)=\bar t$.
\end{cpf}

To complete the proof we need to establish that $\hat \pi(x) + \bar \pi(-x) \geq 1$. Let $\alpha = \hat \pi(x)$. By Claim~\ref{claim:alpha} and the fact that $-x$ viewed as a real number is $1-x$, it suffices to show that $\sup\{ t: \mu(\pi^{-1}([0,t))) \leq 1-x\} \geq 1 - \alpha$. By property~\ref{two-prop:3} in Lemma~\ref{lem:two-prop}, $\mu(\pi^{-1}([0,\alpha])) \geq x$. By symmetry of $\pi$, this implies that $\mu(\pi^{-1}([1-\alpha,1])) \geq x$, and therefore $1-\mu(\pi^{-1}([0,1-\alpha))) \geq x$. Thus, $\sup\{ t: \mu(\pi^{-1}([0,t))) \leq 1-x\} \geq 1-\alpha$.\endproof
\medskip

We now state our main rearrangement theorem. In its proof, we shall need the following technical lemma about monotone, i.e., nondecreasing or nonincreasing, functions; see~\cite[Theorem 4.30]{babyrudin}.

\begin{lemma}\label{lem:countable-discont} Any monotone real valued function defined on any real interval has countably many discontinuities. 
\end{lemma}

\begin{theorem}\label{thm:final-construction} Let $G$ be any compact, connected topological group with Haar measure $\mu$, and $\bve \in G\setminus\{\0\}$. Let $\pi \in \mathcal M_\bve(G)$. Define $\hat \pi$ as in~\eqref{eq:rearrange} and $\bar \pi$ as in~\eqref{eq:right-limit-fn}. Then the function
\begin{equation}\label{eq:final-rearrange}
\tilde \pi := \frac{\hat \pi + \bar \pi}{2}\end{equation}
defined on $\T^1$ is nondecreasing, subadditive and symmetric in the sense that $\tilde \pi(x) + \tilde \pi(-x) = 1$ for all $x \in \T^1\setminus\{0\}$. Further, $\int_G \ln(\pi)d\mu$ exists and is finite if and only if $\int_{\T^1} \ln(\tilde \pi)d\ell$ exists and is finite, and in this case $$\int_G \ln(\pi)d\mu = \int_{\T^1} \ln(\tilde \pi)d\ell,$$ where $\ell$ denotes the Haar measure on $\T^1$.
\end{theorem}

\proof{Proof.}
Since $\hat \pi$ is nondecreasing by property~\ref{two-prop:2} in Lemma~\ref{lem:two-prop}, and $\bar \pi$ is also nondecreasing by Lemma~\ref{lem:right-limit-fn}, so is $\tilde \pi$. Since $\hat \pi$ and $\bar \pi$ are both subadditive by Theorem~\ref{thm:rearrange} and Lemma~\ref{lem:right-limit-fn}, and subadditivity is preserved by convex combinations, $\tilde \pi$ is subadditive. We now check symmetry of $\bar \pi$:
$$\begin{array}{rclc}
\tilde \pi(x) + \tilde \pi(-x) & = &\frac{\hat \pi(x) + \bar \pi(x) + \hat \pi(-x) + \bar \pi(-x)}{2} & \\
& = & \frac{(\hat \pi(x) + \bar \pi(-x)) + (\hat \pi(-x) + \bar \pi(x))}{2} & \\
& = & \frac{1+1}{2} & \qquad\textrm{by Lemma~\ref{lem:right-limit-fn}} \\
& = & 1&
\end{array}
$$

Since $\hat \pi$ is nondecreasing by property~\ref{two-prop:2} in Lemma~\ref{lem:two-prop},
it has countably many discontinuities by Lemma~\ref{lem:countable-discont}. Therefore, $\bar \pi$ differs from $\hat \pi$ only on a countable set, and the same is true for $\tilde \pi$. Thus, all three functions have the same value of the integral on the torus, if the integral exists. Proposition~\ref{prop:preserve-int} then gives the final conclusion.\endproof
\medskip

\begin{theorem}\label{thm:hard-work}
Let $\mathcal{G}$ be the set of functions $h: \T^1 \to \R_+$ with $h(0)=0$ that are nondecreasing, subadditive, and symmetric in the sense that $h(x) + h(-x) = 1$ for all $x \in \T^1\setminus\{0\}$. Then $-\ln(h)$ is integrable for all $h \in \mathcal{G}$ and 
$$\inf\left\{\int_{\T^1} \ln (h)d\ell : h \in \mathcal{G}\right\} = -1.$$
Moreover, the infimum is attained by the function $g(x) = x$, where the right hand side is interpreted as a real number. \end{theorem}

\proof{Proof.}
For any function $h$ on the torus, we will consider the associated function on $[0,1]$ which is identical to $h$ on $[0,1)$ and takes value $1$ at $1$. Now the integral of such a function on $[0,1]$ is equal to the integral of $h$ on the torus. Thus, without further comment, below we will consider all functions on the torus as functions defined on $[0,1]$. With a slight abuse of notation, if $h$ is subadditive on the torus, we will also say the corresponding function on $[0,1]$ is subadditive.

Let $h: \T^1 \to \R_+$ be a function with $h(0)=0$ that is nondecreasing, subadditive, and symmetric. For every prime number $q$, define a function $\pi^q:\Z/q\Z:\to\R_+$ by setting
\[\pi^q(x)=h\left(\frac x{q-1}\right)\quad \mbox{for $x\in\Z/q\Z$},\]
where $x$ is interpreted as a real number in the fraction $\frac x{q-1}$.

It is clear that $\pi^q(0)=0$.
Moreover, for every $x\in\Z/q\Z$ we have
\[\pi^q(x)+\pi^q(q-1-x)=h\left(\frac x{q-1}\right)+h\left(\frac {q-1-x}{q-1}\right)=1,\]
where we used the fact that $h$ is symmetric. Thus $\pi^q$ is symmetric with respect to $b=q-1$.

We show that $\pi^q$ is also subadditive. If $x,y\in\{0,\dots,q-1\}$ are such that $x+y\le q-1$, then, by subadditivity of $h$,
\[\pi^q(x)+\pi^q(y)=h\left(\frac x{q-1}\right)+h\left(\frac y{q-1}\right)\ge h\left(\frac {x+y}{q-1}\right)=\pi^q(x+y).\]
If, on the contrary, $x+y\ge q$, then
\[\pi^q(x)+\pi^q(y)=h\left(\frac x{q-1}\right)+h\left(\frac y{q-1}\right)\ge h\left(\frac {x+y-q+1}{q-1}\right)\ge h\left(\frac {x+y-q}{q-1}\right)=\pi^q(x+y),\]
where in the last inequality we used the fact that $h$ is nondecreasing.

By Lemma \ref{lem:countable-discont}, $h$ has countably many discontinuities. Then, by definition of Riemann integral,
\[\begin{split}
\int_{\T^1}\ln(h)d\ell &=\int_0^1\ln(h(x))dx\\
&= \lim_{q\to\infty}\frac1{q-1}\sum_{x=1}^{q-1}\ln(\pi^q(x))\\
&= \lim_{q\to\infty}\frac1{q-1}\ln\left(\prod_{x=1}^{q-1}\pi^q(x)\right)\\
&\ge \lim_{q\to\infty}\frac1{q-1}\ln\left(\frac{(q-1)!}{(q-1)^{q-1}}\right)\\
&= \lim_{q\to\infty}\frac1{q-1}\sum_{x=1}^{q-1}\ln\left(\frac x{q-1}\right)\\
&= \int_0^1\ln(x)dx=-1,
\end{split}\]
where in the inequality we used Theorem~\ref{thm:finite-volume}.
\endproof
\medskip

\proof{Proof of Theorem~\ref{thm:infinite-volume}.} For any $\pi \in \mathcal M_\bve(G)$, we consider $\tilde \pi: \T^1 \to \R_+$ as defined in~\eqref{eq:final-rearrange} which, by Theorem~\ref{thm:final-construction}, is nondecreasing, subadditive, symmetric and satisfies
$$\int_G \ln(\pi)d\mu = \int_{\T^1} \ln(\tilde \pi)d\ell,$$
where $\ell$ denotes the Haar probability measure on $\T^1$, if either integral exists and is finite. By Theorem~\ref{thm:hard-work}, $-\ln(\tilde \pi)$ is integrable and $\int_{\T^1} \ln(\tilde \pi)d\ell \geq -1$. As a consequence, $-\ln (\pi)$ is integrable and $\int_{G} \ln (\pi)d\mu \geq -1$ for all $\pi \in \mathcal M_\bve(G)$.  The statement about $\GMI^n_{\bve_i}$ follows from the fact that the integral of the logarithm of $\GMI^n_{\bve_i}$ is equal to the integral of the logarithm of $\GMI^1_{\bve_i}$, which is easily verified to be $-1$: $$ \int_{0}^{\bve_i} \ln\big(\frac{x}{\bve_i}\big)dx + \int_{\bve_i}^1 \ln\big(\frac{1-x}{1- \bve_i}\big)dx = \bve_i\int_{0}^{1} \ln(y)dy + (1-\bve_i)\int_{0}^1 \ln(y)dy = \int_{0}^1 \ln(y)dy = -1$$
\endproof
\medskip

We close this section by observing that the Gomory function $\GMI^n_{\bve_i}$ is not the unique minimizer of $\int_{\T^n}\ln(\pi)d\mu$. Indeed, already for $n=1$, the function defined by $\pi(x):=\GMI^1_{b}(kx)$ for $x\in\T^1$ achieves value of the integral equal to $-1$, for every $k\in\N$.

%
%
%

\section*{Acknowledgments.} We gratefully acknowledge deep insights from an anonymous referee. These helped us to considerably simplify the proofs of Theorems~\ref{thm:finite-volume} and~\ref{thm:infinite-volume} from an earlier version of the paper. We are also thankful to Santanu Dey for pointing us to the papers~\cite{bowman1970finiteness} and~\cite{rubin1972strengthened}.

Amitabh Basu was supported by NSF grant CMMI1452820. Michele Conforti and Marco Di Summa were supported by the grants ``SID 2016'' and ``PRIN 2016''.

\bibliographystyle{plain}
\bibliography{full-bib}

\end{document}